\documentclass[10pt]{amsart}
\usepackage[latin1]{inputenc}
\usepackage{amsmath}
\usepackage{amssymb}
\usepackage{amsfonts}
\usepackage{graphicx}
\usepackage[usenames,dvipsnames]{color}
\usepackage{color}

\def \C {\mathbb C}

\newcommand{\sphere}{\overline{\mathbb{C}}}
\newcommand{\pemlam}{\{ P_{\lambda ;m} \}_{m=1}^{\infty}}
\newcommand{\pemlamo}{\{ P_{\lambda_0 ;m} \}_{m=1}^{\infty}}
\newcommand{\pemlamt}{\{ \tilde P_{\lambda ;m} \}_{m=1}^{\infty}}
\newcommand{\linf}{l^{\infty}(\mathbb{C})}
\newcommand{\phimn}{\phi_{\lambda ;m}^n(z)}
\newcommand{\basinm}{\mathcal{A}_{\infty,m}}
\newcommand{\imag}{{\mathbf i}}

\def \Pm {\{P_m \}_{m=1}^\infty}

\def \Pmt {\{\tilde P_m \}_{m=1}^\infty}

\newcommand{\D}{{\mathbb D}}

\theoremstyle{plain}
\newtheorem{theorem}{Theorem}[section]
\newtheorem{proposition}{Proposition}[section]

\newtheorem*{definition*}{Definition}
\newtheorem{lemma}{Lemma}[section]
\newtheorem{corollary}{Corollary}[section]
\newtheorem{corollary*}{Corollary}

\begin{document}
\author{Mark Comerford}
\email{mcomerford@math.uri.edu \vspace{-.3cm}}

\author{Todd Woodard}

\email{twoodard@math.uri.edu \vspace{-.2cm}}

\address{Dept. of Mathematics, University of Rhode Island, Kingston, RI 02881}

\title{Preservation of External Rays in non-Autonomous Iteration
}

\maketitle
\section*{Abstract}

We consider the dynamics arising from the iteration of an arbitrary sequence of polynomials with uniformly bounded degrees and coefficients and show that, as parameters vary within a single hyperbolic component in parameter space, certain properties of the corresponding Julia sets are preserved.  In particular, we show that if the sequence is hyperbolic and all the Julia sets are connected, then the whole basin at infinity moves holomorphically. This extends also to the landing points of external rays and the resultant holomorphic motion of the Julia sets coincides with that obtained earlier in \cite{Com3} using grand orbits. In addition, we have combinatorial rigidity in the sense that if a finite set of external rays separate the Julia set for a particular parameter value, then the rays with the same external angles separate the Julia set for every parameter in the same hyperbolic component.\footnote{To appear in the Journal of Difference Equations and Applications}

\section{Introduction}

In the theory of iteration of complex polynomials, holomorphic motions and external rays are two of the most important tools used to describe parameter space. Ma\~n\'e, Sad and Sullivan showed in \cite{MSS} that the Julia set moves holomorphically for a dense set of parameters and this work has been partly extended to non-autonomous iteration in \cite{Com3}. Douady and Hubbard showed in \cite{DH} that the Mandelbrot set can be largely described in terms of the landings of external rays, and there has been a great deal of work done recently toward applying this work to connectedness loci for other contexts (see for example \cite{BS, Sum4} and particularly the work of Sester on fibered quadratic polynomials in \cite{Ses2}).   

We show that we may always define a B\"ottcher isomorphism on the complement of the filled Julia set.  If we further assume hyperbolicity and connectedness of the Julia sets (which is only natural in view of the classical theory), we will see that this sequence of convenient conformal isomorphisms of the basins of attraction of infinity behave nicely under perturbation, In fact, we 
show that the whole basin at infinity moves holomorphically in the appropriate sense, as do external rays (Corollary 3.1) and that these motions extend to the Julia sets via the landing points of these rays (Theorem 4.1). Further, these motions of the Julia sets must coincide with those given in \cite{Com3} which were constructed using grand orbits rather than external rays.  Finally, given a point in one of the Julia sets for such a sequence, the external rays landing on that point will have the same angles as those landing on an appropriately defined conjugate point for any sequence with parameters in the same hyperbolic component (see Figures 1 and 2).

\vspace{.3cm}
\begin{figure}
\begin{center}
\includegraphics[width=11.1cm]{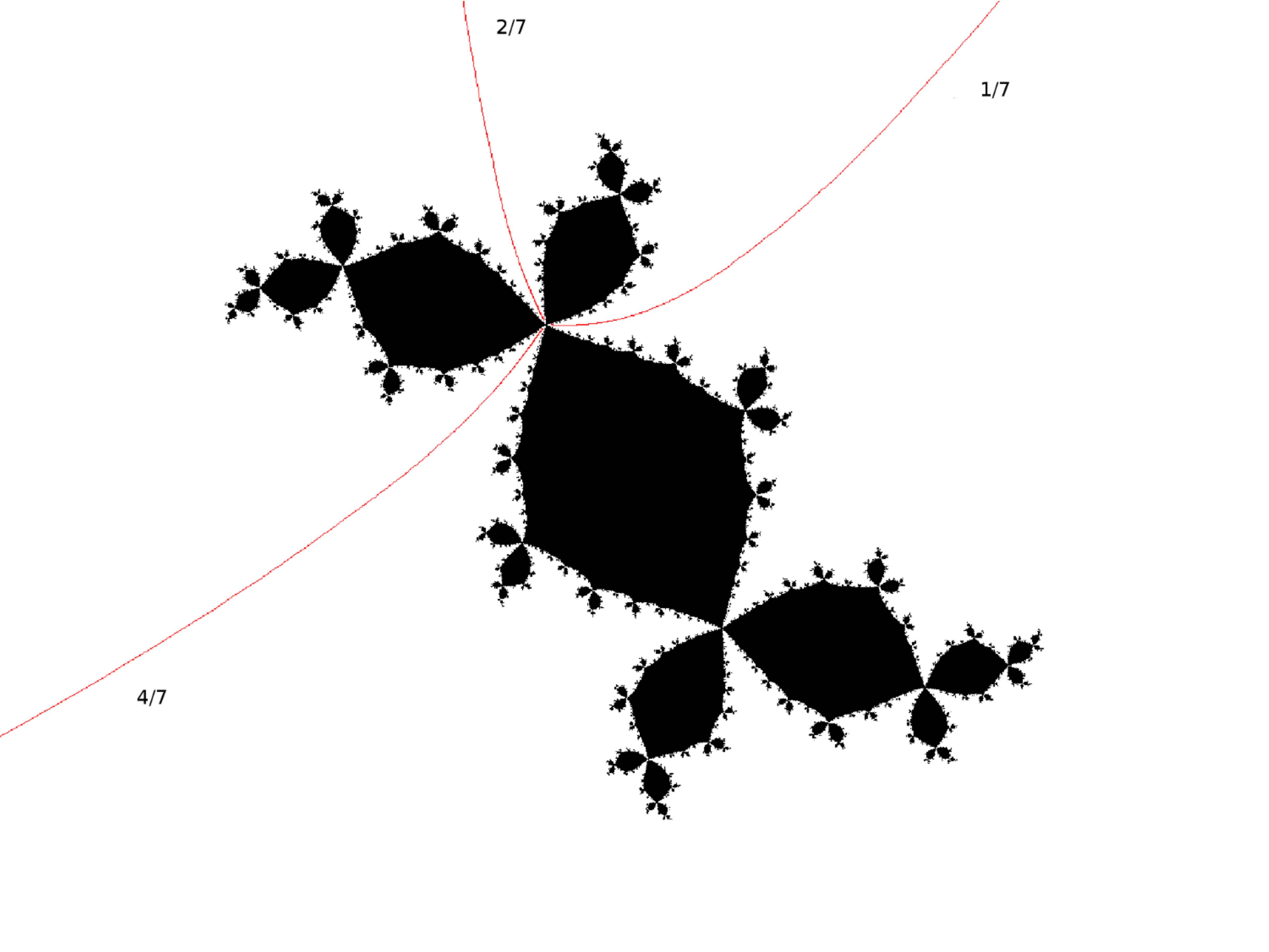}\caption{Classical Julia set for the Rabbit Map $p(z)=z^2  - .745 + .123 \imag$ }
\vspace{1.3cm}
\includegraphics[width=11.1cm]{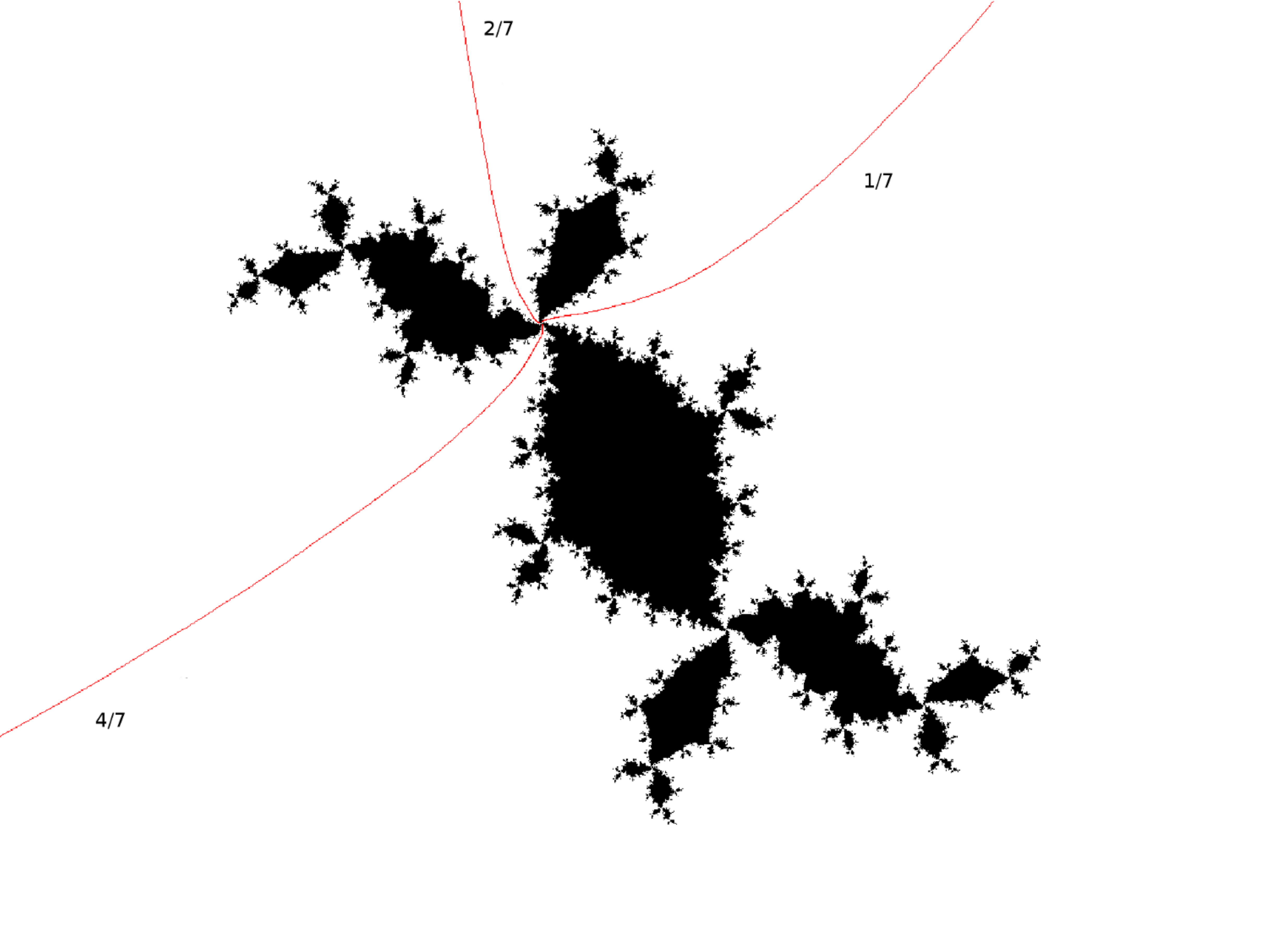}\caption{Non-autonomous Julia set for a sequence of maps $P_m=z^2+c_m$ with $c_m$ selected uniformly at random from a disk of radius 0.06 about $-.745 + .123 \imag$ }
\end{center}
\end{figure}

\vfill
\newpage
This is particularly interesting because of the reliance on periodicity in the classical case and the absence of meaningful periodic orbits in the setting of non-autonomous iteration.  In fact, in the non-autonomous setting, landing points of multiple dynamic rays can frequently be seen as an echo of periodic behavior inherited from a classical (constant) system sharing the same hyperbolic component in parameter space.  Viewing the non-autonomous systems from this perspective gives one hope that characterizing the inherently messy non-autonomous parameter space may not be quite as onerous as it at first seems.

\section{Preliminaries:  Hyperbolic Sequences and Holomorphic Motions}

We will begin with a brief overview of the necessary definitions and basic theory of non-autonomous iteration, as well as two concepts crucial to this paper: hyperbolicity of a non-autonomous sequence of maps, and the holomorphic motion of a set in $\sphere$.

Let $d \ge 2$, $K \ge 1$, $M \ge 0$ and let $\{P_m\}_{m=1}^\infty$ be a sequence of polynomials of the form 

\vspace{-.3cm}
$$
P_m(z)=a_{d_m, m}z^{d_m} +a_{d_m-1, m}z^{d_m-1}+\dots + a_{1, m}z+a_{0, m}
$$
 
such that 
\begin{itemize}
\item $d_m = \deg P_m$ satisfies $2 \le d_m \le d$, for all $m \ge 0$,\\

\item $1/K\leq a_{d_m, m}\le K$ for all $m \ge 0,\\$ 

\item $|a_{n,m}|\le M$ for all $m \ge 0$ and all $0 \le n < d_m$.
     
\end{itemize}

We call such sequences \emph{bounded polynomial sequences} or simply \emph{bounded sequences} and we will refer to the numbers $d$, $K$, $M$ as the \emph{bounds} for the sequence $\Pm$. We will see (in Lemma 3.1) that such a bounded sequence is conformally conjugate (in the appropriate sense, see for example \cite{Com2}) to a polynomial sequence whose members are all monic polynomials.

For $0 \le m$, denote by $Q_m$ the composition $P_m \circ \cdots \cdots \circ P_2 \circ P_1$ and for $0 \le m \le n$ by 
$Q_{m,n}(z)$ the composition $P_{n} \circ \cdots \cdots \circ P_{m+2} \circ P_{m+1}$, (where $Q_{m,m}$ is simply the identity). Let $D_m$ and $D_{m,n}$ denote the degrees of $Q_m$ and $Q_{m,n}$ respectively, so that $D_m = \prod_{i=1}^m d_i$ and $D_{m,n} = \prod_{i=m+1}^n d_i$. If $\Pm$ is a bounded sequence, it is easy to see that we can find $R > 0$ such that for all $m \ge 0$, if $|z| > R$, then $|Q_{m,n}(z)| \to \infty$ as $n \to \infty$. Such a radius is the called \emph{an escape radius} for the sequence $\Pm$. Note that we can find an escape radius $R$ which depends only on the bounds $d$, $K$, $M$ for our sequence and which works for every sequence which satisfies these bounds. We then define the sets
\vspace{.2cm}
\begin{eqnarray*}
\mathcal{K}_m & = & \{z \in \mathbb{C}: \limsup_{n \to \infty}|Q_{m,n}(z)| < \infty\}, \\
\basinm & = & \{z \in \mathbb{C} : \lim_{n \to \infty} |Q_{m,n}(z)| = \infty\},\\
\mathcal{J}_m & = & \partial \mathcal{K}_m = \partial \basinm.
\end{eqnarray*}

where ${\mathcal K}_m$ is called the \emph{$m^{th}$ iterated filled Julia set}, $\mathcal{A}_{\infty,m}$ is the \emph{$m^{th}$ iterated basin of attraction of infinity}, and $\mathcal{J}_m$ is the \emph{$m^{th}$ iterated Julia set}.  As in the classical theory, the \emph{$m^{th}$ iterated Fatou set} $\mathcal{F}_m=\mbox{int}(\mathcal{K}_m) \cup \basinm$ is the domain of normality for the family of functions $\{Q_{m,n}\}_{n=m+1}^\infty.$  It is easy to show that these sets are forward and backward invariant in the sense that for any $0 \le m \le n$, $Q_{m,n}({\mathcal F}_m)={\mathcal F}_n$ and $Q_{m,n}({\mathcal J}_m)={\mathcal J}_n$ and that $Q_{m,n}$ maps components of ${\mathcal F}_m$ onto components of ${\mathcal F}_n$.

We call a bounded sequence of polynomials $\{P_m\}_{m=1}^\infty$ \emph{hyperbolic} if it is uniformly expanding on its iterated Julia sets; that is, if there exist constants $C>0, \mu>1$ such that for all $i,m \ge 0$ and $z \in \mathcal{J}_m$,

\vspace{-.4cm}

$$
|Q^\prime_{m,m+i}(z)|\geq C\mu^i.
$$

For convenience, if $\{P_m\}_{m=1}^\infty$ is bounded and hyperbolic as above, we shall refer to the degree and coefficient bounds $d$, $K$, $M$ and the numbers $C$, $\mu$ as the \emph{hyperbolicity bounds} associated with $\{P_m\}_{m=1}^\infty$.

As in classical complex dynamics, hyperbolicity is an extremely strong condition and, even in this generalized setting, a great deal can be said about sequences with this property.  We quote here a number of results on hyperbolic sequences which will be of use to us in this paper.  The first theorem states that, as in classical complex dynamics, hyperbolicity is an open and stable condition. As each polynomial in a bounded sequence has only finitely many coefficients, it is natural to express the coefficients of the entire sequence of polynomials as a vector $\{\{a_{n,m}\}_{n=0}^{d_m}\}_{m=1}^{\infty}$in $\linf$ where we first list the coefficients of $P_1$, then those of $P_2$ and so on. 

\begin{theorem}  [\cite{Com4} Theorem 3.2] Let $\{P_m\}_{m=1}^ \infty$ be a bounded hyperbolic sequence. Then there exists an open neighbourhood of $\{\{a_{n,m}\}_{n=0}^{d_m}\}_{m=1}^{\infty}$ in $\linf$ such that the bounded sequence arising naturally from any element of this neighbourhood is also hyperbolic.  Moreover, the Julia sets $\mathcal{J}_m$ move continuously (in the Hausdorff topology) as the sequence $\{\{a_{m,n}\}_{n=0}^d\}_{m=1}^{\infty}$ varies throughout this neighbourhood.\end{theorem}

It turns out that one can choose the neighbourhood in the statement so that all polynomial sequences in this neighbourhood are uniformly bounded and are hyperbolic with the same constants.
Related to this theorem is the following useful result. We say that a sequence $\{\{P_m^i\}_{m=1}^\infty\}_{i=0}^\infty$ of bounded sequences converges \emph{pointwise} to another limit sequence $\{P_m\}_{m=1}^\infty$ if, for each $m \ge 0$ and each $0 \le n \le d_m$, the coefficients $a_{n,m}^i$ for $P_m^i$ converge to the corresponding coefficient $a_{n,m}$ of $P_m$.

\vspace{.1cm}
\begin{theorem}  [\cite{Com4} Theorem 3.4] Let $\{\{P_m^i\}_{m=1}^\infty\}_{i=0}^\infty$ be a sequence of bounded polynomial sequences which are bounded and hyperbolic with the same constants and let ${\mathcal J}^i_m$ be the corresponding iterated Julia sets. Suppose also that $\{P_m^i\}_{m=1}^\infty$ converges pointwise to a bounded sequence $\{P_m\}_{m=1}^\infty$ with iterated Julia sets ${\mathcal J}_m$. Then, for each $m \ge 0$, ${\mathcal J}^i_m \to {\mathcal J}_m$ in the Hausdorff topology as $i \to \infty$.
\end{theorem}

\vspace{-.2cm}
This also follows from a result of Sumi (\cite{Sum1} page 583 Theorem 2.14) as well as a result in the paper of Sester (\cite{Ses1} page 411, Proposition 4.1), both of whom were working in the context of polynomials fibered over a compact set. The primary usefulness of this and other similar pointwise limit theorems lie in the construction of Cantor diagonalization arguments like the one in the proof of Lemma 4.1. 

If $\Pm$ is a bounded sequence, for each $0 \le m < n$, let us denote by $C_{m,n}$ the set of critical values of $Q_{m,n}$ which is a set at time $n$. We then define the \emph{postcritical distance}, $PD(\Pm)$, by 

\[ PD(\Pm) = \inf_{m \ge 0, n > m} \mathrm{dist}(C_{m,n} , {\mathcal J}_n)\]

where $dist\,(\cdot,\cdot)$ is the usual Euclidian distance between sets. We will need the following result whose proof can be found in \cite{Com4}. We remind the reader that one condition is said to imply another \emph{up to constants} if the constants associated with the first condition give non-trivial bounds for those associated with the second. 

\vspace{.2cm}
\begin{theorem}  [\cite{Com4} Theorem 1.3] Let $\{P_m\}_{m=1}^{\infty}$ be a bounded sequence.  Then $\{P_m\}_{m=1}^{\infty}$ is hyperbolic if and only if $PD(\Pm) \ge \delta$ for some $\delta >0$. Furthermore, this equivalence is up to constants.
\end{theorem}

This result is the non-autonomous analogue of the result from classical complex dynamics that hyperbolicity is equivalent to the closure of the postcritical set being disjoint from the Julia set and also follows from the work of Sester (\cite{Ses1} page 395 Th\'eor\`eme 1.1). For complete proofs of these statements and a more detailed treatment of the Fatou-Julia theory in the setting of non-autonomous iteration, the reader is referred to \cite{BB, Com1, Com3, Com4}. 

We make extensive use of the theory of holomorphic motions, originally developed by Ma\~ne, Sad, and Sullivan in \cite{MSS}, which, following \cite{Com3}, we consider in the case of analytic dependence on a parameter in $\linf$.  Recall that if $\Lambda \subset \linf$ is open, a function $f: \Lambda \mapsto \mathbb C$ is \emph{G\^ateaux-holomorphic} at $\lambda \in \Lambda$ if for each vector $\xi$ in $\linf$, the complex-valued function of one complex variable $t \longmapsto f(\lambda + t\xi)$ is holomorphic. The linear operator ${\mathbf d}_\xi f$ is called the \emph{differential} of $f$ acting on $\xi$ at the point $\lambda$. If this operator is bounded, and we have that 

\vspace{-.4cm}
\[||f(\lambda + \xi) -f(\lambda) - {\mathbf d}_\xi f(\lambda)|| = o(||\xi||)\]

as $\xi \to 0$, then we say that $f$ is \emph{Fr\'echet-holomorphic} or simply {\it holomorphic} at $\lambda$ and we denote the Fr\'echet differential by ${\mathbf d}_{\lambda}f$. 

We say that $f$ is \emph{Fr\'echet-holomorphic} on $\Lambda$ if it is Fr\'echet-holomorphic at each point of $\Lambda$. This is equivalent to $f$ being G\^ateaux-holomorphic and continuous
on $\Lambda$. Finally, it is well known that the condition of 
continuity can be replaced by that of local boundedness. See \cite{Dineen} for details.

Let $X \subseteq \sphere$ and $\Lambda$ be an open ball in some normed complex vector space $V$ (usually $\ell^\infty(\mathbb{C})$) with centre $\lambda_0$.  We say that a collection of functions $\{\tau_\lambda(x):\lambda \in \Lambda, x \in X\}$ is a  \emph{holomorphic motion of the set $X$} if the following hold:

\begin{itemize}

\item[(i)] $\tau_\lambda(x)$ is continuous and injective as a function of $x$ with $\lambda$ held fixed;\\
\item[(ii)] $\tau_\lambda(x)$ is Fr\'echet-holomorphic as a function of $\lambda$ with $x$ held fixed;\\
\item[(iii)] $\tau_{\lambda_0}(x)=x$ for all $x \in X.$

\end{itemize}

In the case that $V$ is not locally compact (as is the case with $\ell^\infty(\mathbb{C})$), we also have the following technical requirement:

\begin{itemize}

\item[(iv)] There are three distinct points $x_1,x_2,x_3 \in X$ and a positive real number $r$ such that 

$$
\inf_{\lambda \in {\mathrm B}(\lambda_0,r')} \mbox{d}^{\#}(\tau_\lambda(x_i),\tau_\lambda(x_j))>0 \ \ \ \ \ \ \forall \ r', \ 0<r'<r \ \mbox{and} \ i\neq j.
$$  

\end{itemize}

where $d^\#$ denotes the spherical metric on $\sphere$.  The famous $\lambda$-lemma in \cite{MSS} establishes that if $T:\mathbb{D}\times X\rightarrow \sphere$ is a holomorphic motion, then there exists a quasiconformal extension of $T$, $\overline{T}:\mathbb{D}\times \overline{X}\rightarrow \sphere$.  This result allows the authors to prove that classical Julia sets move holomorphically through hyperbolic maps, and hence rational functions arising from the same hyperbolic component are quasiconformally conjugate on their Julia sets.  

The analogous result to the $\lambda$-lemma appropriate to the present context (where $\mathbb{D}$ is replaced by a ball in $\ell^\infty(\mathbb{C})$) is given in \cite{Com3}, Theorem 2.2. If for each $m \ge 0$ we have subsets $U_m$, $V_m$ of $\sphere$ and a continuous mapping $\phi_m: U_m \mapsto V_m$, then we say that the sequence $\{\phi_m\}_{m=0}^\infty$ is \emph{equicontinuous} on $\{U_m\}_{m=0}^\infty$ if for every $\epsilon > 0$, there
exists $\delta >0$ such that for all $m \ge 0$ and all
$x_m, y_m$ in $U_m$ with ${\rm d}^\#(x_m, y_m) < \delta$, we have ${\rm d}^\#(\phi_m(x_m), \phi_m(y_m)) < \epsilon$
where ${\rm d}^\#(.\,,.)$ refers to the standard spherical distance between points in $\sphere$. We say $\{\phi_m\}_{m=0}^\infty$ is \emph{bi-equicontinuous} on $\{U_m\}_{m=0}^\infty$ if in addition each $\phi_m$ is a bijection from $U_m$ to $V_m$ and the family of inverse mappings $\{\phi_m^{\circ -1}\}_{m=0}^\infty$ is equicontinuous on $\{V_m\}_{m=0}^\infty$.

In order to utilize these ideas, we require the concept of an \emph{analytic family} of sequences of polynomials depending on an infinite-dimensional complex parameter. This definition is taken from \cite{Com3}.

\begin{definition*}Let $\Lambda \subset l^\infty$ be open.
We say a family of bounded sequences\\ ${\mathcal P}_{\Lambda} = \{ \pemlam, \lambda \in \Lambda\}$ is an {\rm analytic family over $\Lambda$} if:

\vspace{-.2cm}
\begin{enumerate}
\item For each $\lambda$, $P_{\lambda;m}$ is a bounded polynomial sequence and these bounds are locally uniform
on $\Lambda$;

\vspace{.2cm}
\item The coefficients $a_{i,m}(\lambda)$ 
are all Fr\'echet-differentiable and their Fr\'echet differentials are locally uniformly bounded on $\Lambda$.
\end{enumerate}

\end{definition*}

For convenience, if ${\mathcal P}_{\Lambda}$ is an analytic family of bounded sequences as above, if the degree and coefficient bounds $d$, $K$, $M$ are uniform, as well as the bound $B$ on the operator norms of the Fr\'echet differentials of the coefficients of all polynomials involved, we shall refer to $d$, $K$, $M$, $B$ as the \emph{constants} associated with this analytic family. Also, if $\Lambda_0 \subset \Lambda$, then we shall call the family $\{ \pemlam, \lambda \in \Lambda_0\}$ 
the \emph{restriction of ${\mathcal P}_{\Lambda}$ to $\Lambda_0$} and denote it by 
${\mathcal P}_{\Lambda_0}$. Finally, if all the polynomials of all the sequences of an analytic family ${\mathcal P}_{\Lambda}$ as above are monic, then we will call ${\mathcal P}_{\Lambda}$ a \emph{monic analytic family}.

We call the set of all points $\lambda \in \Lambda$ for which all the iterated Julia sets 
${\mathcal J}_{\lambda; m}$, $m \ge 0$ are connected, the \emph{connectedness locus} for ${\mathcal P}_{\Lambda}$ and denote it by ${\mathcal M}_{\Lambda}$. We denote the set of all $\lambda$ for which the iterated Julia sets ${\mathcal J}_{\lambda; m}$ are connected and the sequence $\pemlam$ is hyperbolic by $\mathcal{HC}_{\Lambda}$. Lastly, for $\lambda \in \linf$, $\epsilon > 0$, we denote the ball $\{\lambda:||\lambda-\lambda_0||_{\infty} < \epsilon \}$ by ${\mathrm B}(\lambda_0, \epsilon_0)$.

\begin{theorem}[\cite{Com3}, Theorem 1.3] Let $\Lambda \subset \linf$ be open
and let ${\mathcal P}_{\Lambda}$ be an analytic family of bounded polynomial sequences over $\Lambda$. Let $\lambda_0 \in \Lambda$ and suppose that $\{P_{\lambda_0;m}\}_{m=1}^\infty$ is bounded and hyperbolic. Suppose also that the operator norms of the derivatives of the coefficients are bounded by $B$ on an $\epsilon$-neighbourhood of $\lambda_0$ for some $\epsilon >0$. Then there exist neighbourhoods $\Lambda_0$ in $\linf$ of $\lambda_0$ in $\Lambda$, and for each $m \ge 0$ 
an embedding $\tau_m(\lambda, z) = (\lambda, \tau_{\lambda;m}(z))$ defined for $\lambda \in \Lambda_0$ and $z \in {\mathcal J}_{\lambda_0; m}$ for which we have the following:
\vspace{-.2cm}
\begin{enumerate}
\item $\Lambda_0 = {\mathrm B}(\lambda_0, \epsilon_0)$ is an $\epsilon_0$-neighbourhood of $\lambda_0$ where $\epsilon_0$ depends on $\epsilon, B$, and the hyperbolicity bounds for $\Pm$;

\vspace{.25cm}
\item $\tau_{\lambda;m}$ is holomorphic in $\lambda$ and the Fr\'echet derivative ${\bf d}\tau_{\lambda;m}$ has operator norm bounded uniformly on $\Lambda_0$, the bounds depending only on $B$ and the hyperbolicity bounds for $\Pm$ (and in particular not on $m$, $z$, or $\lambda$);

\vspace{.25cm}
\item $\tau_{\lambda;m}$ is quasiconformal and bi-equicontinuous in $z$ on ${\mathcal J}_{\lambda_0; m}$, and jointly continuous as a function of $z$ and $\lambda$ on ${\mathcal J}_{\lambda_0; m} \times \Lambda_0$, the estimates depending only on $\frac{||\lambda - \lambda_0||}{\epsilon_0}$, $B$, and 
the hyperbolicity bounds for $\Pm$;  

\vspace{.25cm}
\item $\tau_{\lambda;m}({\mathcal J}_{\lambda_0; m}) = {\mathcal J}_{\lambda; m}$ for every 
$\lambda \in \Lambda_0$ and each $m \ge 0$;  

\vspace{.25cm}
\item $\tau_{\lambda_0;m} = {\rm Id}_{{\mathcal J}_{\lambda_0; m}}$ and $\tau_{\lambda, m+1}\circ P_{\lambda_0;m} = P_{\lambda;m} \circ \tau_{\lambda;m}$ on 
${\mathcal J}_{\lambda_0; m}$ for each $m \ge 0$.
\end{enumerate}
 
\end{theorem} 

The holomorphic motion above was constructed using grand orbits which gave a dense subset of the iterated Julia sets on which one could then apply the $\lambda$-Lemma. However, the result is local in nature and the holomorphic motion could also potentially depend how it was constructed and, even in the case of grand orbits as above, 
on the choice of grand orbit. The next results show that this is not the case. We remind the reader that, since, by Theorem 2.1, hyperbolicity is an open condition, it makes sense to speak of hyperbolic components for an analytic family over some open subset $\Lambda$ of $\linf$.

\vspace{.2cm}
\begin{theorem}
Let $\Lambda \subset \linf$ be open
and let ${\mathcal P}_{\Lambda}$ be an analytic family of bounded polynomial sequences over $\Lambda$. Let $\lambda_0 \in \Lambda$
and let $\Lambda_0 = {\mathrm B}(\lambda_0, \epsilon) \subset \Lambda$ be such that all the sequences $\pemlam$, $\lambda \in \Lambda_0$ are bounded and hyperbolic with the same bounds. Let $\sigma_{\lambda; m}(z)$, $\tau_{\lambda; m}(z)$ be 
two holomorphic motions  defined on the iterated Julia sets ${\mathcal J}_{\lambda; m}$, $\lambda \in \Lambda_0$, which are conjugacies on the iterated Julia sets in the sense that (4) and (5) of Theorem 2.4 above hold. 

Then there exists $0 < \epsilon'  \le \epsilon$ so that, for all $\lambda \in \Lambda_1 := {\mathrm B}(\lambda_0, \epsilon')$, $m \ge 0$ and $z \in {\mathcal J}_{\lambda_0; m}$, $\sigma_{\lambda; m}(z) = \tau_{\lambda; m}(z)$. 
\end{theorem}

\textbf{Proof:} Since all the sequences $\pemlam$ are bounded with the same bounds, we may find a uniform escape radius $R_0$ for all all them. By setting $\epsilon' = \epsilon/2$, and applying the Schwarz lemma, we see that for any $\lambda \in {\mathrm B}(\lambda_0, \epsilon')$ and any fixed $z \in {\mathcal J}_{\lambda; m}$, the G\^ateaux derivatives of $\sigma_{\lambda; m}(z)$, $\tau_{\lambda; m}(z)$ in any fixed direction given by a unit vector $\xi$ in $\linf$ are uniformly bounded above by $2R_0/\epsilon$. Hence the operator norms of the Fr\'echet differentials ${\mathbf d}_\lambda \sigma_{\lambda; m}(z)$, 
${\mathbf d}_\lambda \tau_{\lambda; m}(z)$ are bounded on $\Lambda_1 = {\mathrm B}(\lambda_0, \epsilon')$ and that these bounds are uniform with respect to the bounds for the restricted family ${\mathcal P}_{\Lambda_1}$ and do not depend on $m$, $z$ or $\lambda$. 

By Theorem 2.3 and the uniform hyperbolicity of the statement, we may find $\delta > 0$ such that for every $\lambda \in \Lambda_0$, the sequence $\pemlam$ has postcritical distance $\ge \delta$. Again by uniform hyperbolicity, we may find $N_0 \ge 1$ such that for every choice of $\lambda \in \Lambda_0$, $m \ge 0$, and $z \in {\mathcal J}_{\lambda; m}$, we have $|Q_{\lambda; m, m + N_0}'(z)| \ge 2$ (such a number is called a \emph{doubling time} for our sequences, see e.g. \cite{Com4}). It then follows easily from the standard distortion theorems for univalent mappings (e.g. \cite{CG} P. 3, Theorem 1.6) that we may find a universal constant $0 < c < 1$ such that for every such $\lambda$, $m$, $z$, we have the following two things:

\vspace{-.1cm}
\begin{enumerate}

\item $Q_{\lambda; m, m + N_0}$ is univalent on the disc ${\mathrm D}(z, c\delta)$;

\vspace{.3cm}
\item The component of the inverse image $Q_{\lambda; m, m + N_0}^{\circ -1} ({\mathrm D}(Q_{\lambda; m, m + N_0}(z), c\delta))$ which contains $z$ is contained in ${\mathrm D}(z, c\delta)$.
\end{enumerate} 

\vspace{-.1cm}
Now fix $m \ge 0$ and $z \in {\mathcal J}_{\lambda_0; m}$. For each $n \ge m$ let $\Delta_{\lambda; n} = {\mathrm D}(\sigma_{\lambda; n}(z), c\delta)$ be a ball of radius $c\delta$ about $\sigma_{\lambda; n}(z)$. 

Our two holomorphic motions must agree at $\lambda_0$ and the argument is based on the fact that they do not move apart too `fast' in $\lambda$ and one can then use the contraction of inverse branches of the iterates which comes from hyperbolicity to show that they do not separate at all. By the uniform bounds on the Fr\'echet differentials of $\sigma_{\lambda; m}(z)$, $\tau_{\lambda; m}(z)$, we can make $\epsilon'$ smaller if needed such that for every $\lambda \in {\mathrm B}(\lambda_0, \epsilon')$ and every $n \ge m$ we have $\tau_{\lambda; n}(z) \in \Delta_{\lambda; n}$. 

Recall the degrees $D_{m,n}={\prod_{i=m+1}^n d_i}$ of the compositions $Q_{\lambda; m,n}$ (note that, as observed in \cite{Com3}, by the definition of an analytic family, it is easy to see that these will not depend on $\lambda$). The inverse image of $\Delta_{\lambda; m + N_0}$ must consist of $D_{m, m+N_0}$ distinct components, only one of which contains $\sigma_{\lambda; m}(z)$. This component also contains a preimage of $\tau_{\lambda; m + N_0}(Q_{\lambda_0; m, m + N_0}(z))$ and it follows from our two requirements for $c$ above that this preimage must agree with $\tau_{\lambda; m}(z)$. 

A simple induction involving repetition of this argument shows that, for any  $k \ge 1$, the preimage of $\tau_{\lambda; m + kN_0}(Q_{\lambda_0; m, m + kN_0}(z))$ in the component of $Q_{\lambda; m , m+ kN_0}^{\circ -1}(\Delta_{\lambda; m + kN_0})$ which contains $\sigma_{\lambda; m}(z)$ must also agree with $\tau_{\lambda; m}(z)$. It then follows again by the distortion theorems that we can find $C' > 0$ such that $|\sigma_{\lambda; m}(z) - \tau_{\lambda; m}(z)| < C' 2^{-k}$. Letting $k$ tend to infinity then gives the result. $\blacksquare$

\section{B\"ottcher Maps}

In classical complex dynamics much of the geometry of the connectedness locus for the Julia sets in quadratic parameter space is understood primarily in terms of external rays and parameter wakes.  The development of the theory in the classical setting, however, relies heavily on periodic orbits, and specifically parabolic dynamics.  In \cite{Mi2}, it is established that each hyperbolic component of the quadratic connectedness locus $\mathcal{M}$ can be put into bijective correspondence with a parabolic \emph{root point}.  This root point is the landing point of exactly two parameter rays which are periodic under the doubling map, and the period of these rays is a multiple of the period of the unique attracting orbit associated with the corresponding orbit portrait.  

Despite the lack of periodicity in our current setting, we aim to bring some of these same tools to bear in the non-autonomous context, namely B\"ottcher isomorphisms and Green's functions with pole at infinity on the iterated basins of attraction of infinity.  Some work in this direction has already been done.  For example, in \cite{Com4} it is shown that, for each $m \ge 0$, the Green's functions defined by $G_m(z) = \lim_{n \to \infty} (1/D_{m,n}) \log |Q_{m,n}(z)|$ (as originally defined for monic bounded sequences by Fornaess and Sibony in \cite{FS}) exist on the iterated basin of infinity $\mathcal{A}_{\infty, m}$. Moreover, they are continuous with respect to variation of the polynomial sequence by a parameter (\cite{FS}) and satisfy the functional relations $G_n (Q_{m,n}(z)) = D_{m,n}G_m(z)$ on $\mathcal{A}_{\infty, m}$ for every $0 \le m \le n$. 
Finally, Br\"uck shows in \cite{Br} that, for certain bounded sequences of quadratic polynomials, we have a B\"ottcher isomorphism $\phi_m(z)$ between $\mathcal{A}_{\infty, m}$ and $\overline {\mathbb C}\setminus \overline {\mathbb D}$ for which $\phi_{m+1}\circ P_{m+1} \circ \phi_m^{\circ -1} = z^2$ on $\mathcal{A}_{\infty, m}$ for each $m \ge 0$.

Before we prove the main theorem of this paper, we require the following lemma which will allow us to simplify our subsequent calculations somewhat. Recall from \cite{Com2} that two bounded polynomial sequences $\Pm$, $\Pmt$ are said to be \emph{analytically conjugate on $\sphere$} if we can find $A \ge 1$, $B \ge 0$ and a sequence of affine mappings $\{\chi_m\}_{m=0}^\infty$, 
where $\chi_m = \alpha_m z + \beta_m$ satisfies $1/A \le |\alpha_m| \le A$, $|\beta_m| \le B$ for every $m \ge 0$. Note that such a sequence of mappings will be bi-equicontinuous on $\sphere$.

\begin{lemma}
Let $\Lambda \subset \linf$ be open, let ${\mathcal P}_\Lambda$ be an analytic family of bounded polynomial sequences over $\Lambda$ and let $\lambda_0 \in \Lambda$. Then we can find $\epsilon_0 > 0$, a neighbourhood $\Lambda_0 = {\mathrm B}(\lambda_0, \epsilon_0) \subset \Lambda$, a monic analytic family $\tilde {\mathcal P}_\Lambda$ of over $\Lambda_0$ and a sequence $\{\chi_{\lambda; m}\}_{m=0}^\infty$ of affine linear maps for which we have the following:

\vspace{-.2cm}
\begin{enumerate}

\item There exist $A \ge 1$, $C \ge 0$ such that for every $\lambda \in \lambda_0$ and each $m \ge 0$, $\chi_{\lambda;m} = \alpha_m(\lambda) z$ satisfies 

\vspace{.2cm}
\begin{enumerate}
\item $1/A \le |\alpha_m(\lambda)| \le A$;
\vspace{.25cm}
\item The Fr\'echet differential ${\mathbf d}_\lambda \alpha_m (\lambda)$ satisfies
$||{\mathbf d}_\lambda \alpha_m (\lambda)|| \le C$;

\end{enumerate}

\vspace{.25cm}
\item For each $\lambda \in \Lambda_0$, $\pemlam$ is analytically conjugate to $\pemlamt$ via $\{\chi_{\lambda; m}\}_{m=0}^\infty$.

\end{enumerate}

\end{lemma}

\textbf{Proof:} By the definition of an analytic family, we may choose $\epsilon_0 >0$ such that we can find $B \ge 0$ such that on ${\mathrm B}(\lambda_0, \epsilon_0)$ every sequence 
has bounds $d$, $K$, $M$, $B$ as above.

As the Fr\'echet derivatives of all the coefficients are uniformly bounded, we may then make $\epsilon_0$ smaller if needed so that for every $m \ge 0$ we can ensure that for all $\lambda \in \Lambda_0$ the leading coefficients $a_{d_{m}, m}(\lambda)$ avoid either the positive imaginary axis or the negative imaginary axis (or both). We may then take an appropriate branch cut and use a branch of the argument in the range $(-\tfrac{3\pi}{2}, \tfrac{\pi}{2})$ or $(-\tfrac{\pi}{2}, \tfrac{3\pi}{2})$ as appropriate to define the roots $\sqrt[d]{a_{d_m, m}(\lambda)}$ for any integer $d \ge 1$. It is important to note that these assignments of roots are `consistent'. To be specific, if for example we use the branch $(-\tfrac{3\pi}{2}, \tfrac{\pi}{2})$ for our argument, so that $a_{d_m, m}(\lambda) = |a_{d_m, m}(\lambda)|e^{i \theta}$, $-\tfrac{3\pi}{2} < \theta < \tfrac{\pi}{2}$, then for any $d \ge 1$, we define $\sqrt[d]{a_{d_m, m}(\lambda)} := \sqrt[d]{|a_{d_m, m}(\lambda)|}e^{i \theta/d}$ where of course $-\tfrac{3\pi}{2d} < \theta < \tfrac{\pi}{2d}$. Note that then, for $d, d' \ge 1$, $\left (\sqrt[(d d')]{a_{d_m, m}(\lambda)}\right )^{d} = \sqrt[d']{a_{d_m, m}(\lambda)}$ and in particular one has $\left (\sqrt[d]{a_{d_m, m}(\lambda)}\right )^d = a_{d_m, m}(\lambda)$.

For each $m$, define $\alpha_m(\lambda)$ by 
\begin{equation}\alpha_m(\lambda) = \prod_{n=m}^\infty \left (a_{d_n,n}(\lambda)\right )^{1/D_{m-1,n}}.  \end{equation}
Note that for each fixed $n \ge 0$ we use the \emph{same} branch of the argument to form the roots of $a_{d_n,n}(\lambda)$ which appear in $\alpha_m$ as above for every $m \le n$. 
Then it is easy to check that for all $\lambda \in \Lambda_0$, these infinite products will converge with $1/K \le |\alpha_m(\lambda)| \le K$ and that the sequence $\chi_{\lambda;m}$ defined in this way does indeed give us a non-autonomous conjugacy between $\pemlam$ and a bounded monic sequence $\pemlamt$.

The functions $\alpha_m(\lambda)$ are uniformly bounded above and below and if, as in the proof of Theorem 2.5, we replace $\epsilon$ with $\epsilon/2$ in the definition of $\Lambda_0$, then, as these products converge uniformly on $\Lambda_0$, these functions are G\^ateaux-differentiable and their G\^ateaux derivatives are uniformly bounded above by $C:=2K/\epsilon$. The functions $\alpha_m(\lambda)$ are thus Fr\'echet-differentiable and the operator norms of their derivatives are bounded on $\Lambda_0$ by $C$. By the upper and lower bounds on the absolute values of the coefficients $\alpha_m(\lambda)$, similar arguments show that $\pemlamt$ does indeed give us an analytic family over $\Lambda_0$ as required. 
$\blacksquare$

We say two analytic families ${\mathcal P}_\Lambda$, $\tilde {\mathcal P}_\Lambda$ over the same set $\Lambda \subset \linf$ are \emph{globally analytically conjugate} if we can find a sequence $\{\chi_{\lambda; m}\}_{m=0}^\infty = \{\alpha_m(\lambda) z + \beta_m(\lambda)\}_{m=0}^\infty$ of affine linear maps for which the absolute values of the coefficients $\alpha_m(\lambda)$ are locally uniformly bounded above and below and the absolute values of the coefficients $\beta_m(\lambda)$ are locally uniformly bounded above on $\Lambda$. Note that by the Schwarz lemma this also implies that the Fr\'echet derivatives of these coefficients are also locally uniformly bounded on $\Lambda$. 

It is important to note that one cannot use Lemma 3.1 to show that a general analytic family is not automatically conjugate to a monic family as this is strictly a local result. However, for families over a simply connected set $\Lambda$ of parameters (e.g. a ball or more generally a convex set), we can show this.

\vspace{.2cm}
\begin{proposition}
Let $\Lambda \subset \linf$ be simply connected and let ${\mathcal P}_\Lambda$ be an analytic family defined on $\Lambda$. The ${\mathcal P}_\Lambda$ is globally analytically conjugate to a monic analytic family $\tilde {\mathcal P}_\Lambda$. 
\end{proposition}

\vspace{-.2cm}
\textbf{Proof:} Let $\lambda_0 \in \Lambda$, let $\lambda_1$ be any other point of $\Lambda$ and let $\gamma : [0,1] \mapsto \Lambda$ be a path in $\Lambda$ from $\lambda_0$ to $\lambda_1$. For each $t \in [0,1]$, if we let $\lambda_t:=\gamma(t)$, then we have a ball $B(\lambda_t, \epsilon_t)$ about $\lambda_t$ on which the bounds for the sequences $\pemlam$ and the Fr\'echet derivatives of the coefficients of all their polynomials are uniform. 

The image of $\gamma$ is compact as it is the continuous image of a compact set and so we may cover it with a union of finitely many such balls which we denote by $X$. 
We then have a neighbourhood of $\gamma$ on which we have uniform bounds for the sequences and the Fr\'echet derivatives of the coefficients. It also follows from compactness that we can find $\epsilon > 0$ such that for any $0 \le t \le 1$, $B(\lambda_t, \epsilon) \subset X$.

Without loss of generality the first of the balls which make up $X$ is centered about $\lambda_0$. By making this ball smaller if needed we may use Lemma 3.1 to define a local analytic conjugacy $\{\chi_{\lambda; m}\}_{m=0}^\infty$ to a monic family on this ball and we can use this same ball for any path from $\lambda_0$ to $\lambda_1$. Using the uniform bound on the Fr\'echet derivatives of the coefficients of all polynomials arising from sequences with parameters in $X$, it is then not too hard to see that by making all these balls smaller if necessary (in which case we may need more balls, but still finitely many of them), we can choose our branches of the logarithm for the leading coefficients appropriately so as to make use of (1) in the proof of the last result to continue this conjugacy analytically so as to define it on each ball of $X$. Moreover, we can ensure that this continuation of the conjugacy will agree on any overlap between successive balls of $X$. 

From above, our conjugacy is locally uniquely defined in a neighbourhood of $\lambda_0$. Thus, by applying the identity principle on one-dimensional complex lines of the form $\lambda_{t_1} + w (\lambda_{t_2} - \lambda_{t_1})$ for nearby values $t_1, t_2$ of $t$ (where $w \in \C$), it follows that for the points $\lambda_t$ of the image of $\gamma$, the conjugacies 
$\{\chi_{\lambda_t; m}\}_{m=0}^\infty$ will not depend on the choice of discs used to make the set $X$.  In order to obtain a well-defined conjugacy on all of $\Lambda$, we still need to check that the conjugacy $\{\chi_{\lambda_1; m}\}_{m=0}^\infty$
does not depend on the choice of path from $\lambda_0$ to $\lambda_1$ and the argument for this is similar to the classical monodromy theorem from complex analysis (e.g. \cite{Lang} page 311, Theorem 1.3). 

Since $\Lambda$ is simply connected, any two paths $\gamma$, $\eta$ from $\lambda_0$ to $\lambda_1$ are fixed endpoint homotopic in $\Lambda$. Let $\Phi[t,u]: [0,1]^2 \mapsto \Lambda$ be such a homotopy where $\Phi(t,0) = \gamma(t)$, $\Phi(t,1) = \eta(t)$. For each $0 \le u \le 1$, if we consider the path $\gamma_u(t)$ given by $t \mapsto \Phi(t,u)$, we can make a corresponding neighbourhood $X_u$ as above on which we can analytically continue our locally defined conjugacy. We can also as above find $\epsilon_u > 0$ such that $X_u$ is an $\epsilon_u$ neighbourhood of $\gamma_u$ and hence by the continuity of $\Phi$ we can find $\delta_u > 0$ such that the image of $\gamma_s$ is a subset of $X_u$ for all $0 \le s \le 1$ with $|s - u| < \delta_u$. The same argument using the identity principle as before shows that for such $s$, the continuation of $\{\chi_{\lambda_1; m}\}_{m=0}^\infty$ arising from $X_u$ will be the same as that arising from $X_s$ for all points on the path $\gamma_s(t)$, $0 \le t \le 1$. The same compactness argument as in the classical monodromy theorem then shows that the continuation of this conjugacy along $\gamma$ must agree with that along $\eta$ at $\lambda_1$ which completes the proof. 
$\blacksquare$

We now state and prove the main theorem of this paper. 

\begin{theorem}
Let $\Lambda \subset \linf $ be open and let ${\mathcal P}_{\Lambda}$ be an analytic family of bounded polynomial sequences over $\Lambda$. Let $\epsilon > 0$ be such that ${\mathrm B}(\lambda_0, \epsilon) \subset \mathcal{M}_\Lambda$ and the operator norms 
of the Fr\'{e}chet differentials of the coefficients of the members of each sequence $\pemlam$, $\lambda \in {\mathrm B}(\lambda_0, \epsilon)$ are bounded above by some constant $B$. Then there exists $0 < \epsilon_0 \le \epsilon$ and a neighbourhood $\Lambda_0 := {\mathrm B}(\lambda_0, \epsilon_0)$ of $\lambda_0$ 
where $\epsilon_0$ depends on $\epsilon$, $B$ and the degree and coefficient bounds associated with the sequence $\{P_{\lambda_0; m}\}_{m=1}^\infty$ such that the following hold

\begin{enumerate}

\item[(i)] For each $\lambda \in \Lambda_0$, there exists a sequence of B\"ottcher maps $\Phi_\lambda =\{ \phi_{\lambda; m} \}_{m=0}^{\infty}$ where for each $m \ge 0$ $\phi_{\lambda;m}$ maps ${\mathcal A}_{\lambda; \infty ,m}$ conformally to $\sphere \setminus \overline \D$ with 

\vspace{-.2cm}
$$
\phi_{m+1}\circ P_{m+1} \circ \phi_m^{\circ -1} = z^{d_{m+1}}, \qquad z \in {\mathcal A}_{\lambda; \infty, m}.\vspace{.2cm}
$$
Further, for each $\lambda \in \Lambda_0$, each $m \ge 0$ and each $z \in {\mathcal A}_{\lambda; \infty, m}$, these mappings are Fr\'{e}chet-holomorphic in $\lambda$ with $z$ held fixed. Lastly, if ${\mathcal P}_{\Lambda}$ is a monic analytic family and we require that our B\"ottcher mappings have derivative $1$ at infinity, then these mappings are uniquely defined.

\vspace{.3cm}
\item[(ii)]  If we let $\Psi_\lambda =\{\psi_{\lambda;m}\}_{m=0}^\infty$ denote the corresponding sequence of inverse B\"ottcher mappings, then for each $\lambda \in \Lambda_0$, each $m \ge 0$ and each $z \in \sphere \setminus \overline \D$, $\psi_{\lambda;m}$ is also holomorphic in the $\lambda$ variable with $z$ held fixed. Again, if ${\mathcal P}_{\Lambda}$ is a monic analytic family and we require that our B\"ottcher mappings have derivative $1$ at infinity, then these inverse B\"ottcher mappings are uniquely defined. \\

\item[(iii)] For each $\lambda_1 \in \Lambda_0$ and each $m \ge 0$, if $Z$ is any relatively compact open subset (in $\sphere$ with respect to the spherical topology) of $\mathcal{A}_{\lambda_1;\infty,m}$, then there is an open neighbourhood $\Lambda_1
= {\mathrm B}(\lambda_1, r)$ of $\lambda_1$ 
where $r$ depends on $\epsilon$, $||\lambda_1 - \lambda_0||$, the bounds associated with the restricted family ${\mathcal P}_{\Lambda_0}$ and the minimum value of the Green's function $G_{\lambda_1; m}(z)$ on $Z$
such that $\|\mathbf{d}_\lambda \phi_{\lambda;m}\|$ is bounded on $Z\times \Lambda_1$. This bound is uniform with respect to $\epsilon$, $\tfrac{||\lambda - \lambda_1||}{r}$, the bounds associated with the restricted family ${\mathcal P}_{\Lambda_1}$ and the minimum value of $G_{\lambda_1; m}(z)$ on $Z$.

\vspace{.3cm}

\item[(iv)]
For each $m \ge 0$, $\|\mathbf{d}_\lambda\psi_{\lambda;m}\|$ is bounded on $(\sphere \setminus \overline \D) \times\Lambda_0$. This bound is uniform with respect to $\tfrac{||\lambda - \lambda_0||}{\epsilon_0}$ and the bounds associated with ${\mathcal P}_{\Lambda_0}$. In particular, it does not depend on $m$.

\end{enumerate}

\end{theorem}

\textbf{Proof:}  In view of Lemma 3.1, we can assume without loss of generality that all the sequences $\pemlam$, $\lambda \in {\mathrm B}(\lambda_0, \epsilon)$ are monic. Using the bound $B$ on the Fr\'echet derivatives of the coefficients, by changing the constant $M$ slightly if necessary, we can find $\epsilon_0>0$ depending only on these bounds such that if $\Lambda_0$ is the ball ${\mathrm B}(\lambda_0,\epsilon_0)$ and $\lambda \in \Lambda_0$, then $\pemlam$ is monic and has the same degree bound $d$ and the same bound $M$ on the coefficients of the non-dominant terms. 

Since the bounds on the escape radius $R_0$ for the sequence $\pemlamo$ depend only on $d$ and $M$, we may increase $R_0$ so that it is valid for all $\pemlam$ with $\lambda \in \Lambda_0$.  Also, we make $R_0$ large enough that, for each $\lambda \in \Lambda_0$ and each $m \ge 0$, the leading term in $P_{\lambda; m}$ uniformly dominates the lower order terms;  in particular, we insist that for all $d', 2 \leq d' \le d$ and all $R \geq R_0$,
\vspace{-.4cm}

\begin{equation}R^{d'} \geq 2M \sum_{i=0}^{d'-1}R^i.\end{equation}

Recalling the degrees $D_{m,n}={\prod_{i=m+1}^n d_i}$ of the compositions $Q_{\lambda; m,n}$, from (1), provided $R_0^2 /2 \ge R_0$ we then have the following estimate on the size of the iterates $Q_{\lambda;m,n}(z)$ if $\lambda \in \Lambda_0$ and $|z| = R \ge R_0$:
\vspace{-.2cm}

\begin{equation}\left(\frac{1}{2}\right)^{1 + \sum_{i=m+1}^n D_{i,n}}R^{D_{m,n}} \le |Q_{\lambda; m,n}(z)| \le \left(\frac{3}{2}\right)^{1 + \sum_{i=m+1}^n D_{i,n}}R^{D_{m,n}}.\end{equation}

For each $\lambda \in \Lambda_0$, each $m \ge 0$ and $n \ge m$, define the functions 

\vspace{-.35cm}
$$
\phi_{\lambda ;m}^n(z)=\left(Q_{\lambda; m,n}(z)\right)^{1/D_{m,n}}
$$

on ${\mathcal A}_{\lambda; \infty, m}$, where we take the roots by analytically continuing the standard branch (with branch cut on the negative real axis) from some point on the intersection of the positive real axis with the complement of $\overline {{\mathrm D}(0,R_0)}$ (note that, as $Q_{\lambda; m,n}$ is monic, for $z$ real and very large and positive, the argument of $Q_{\lambda; m, n}(z)$ will be small).  As all the iterated Julia sets are connected, there are no critical points in the iterated basins of infinity and a simple lifting argument then shows that each $\phi_{\lambda ;m}^n$ is holomorphic and univalent as a function of $z$ on ${\mathcal A}_{\lambda; \infty, m}$.

The domains of the $\phimn$ differ in $\lambda$ and so, in order to demonstrate that the B\"ottcher maps vary holomorphically with $\lambda$, we must first establish a locally common domain of definition in order to have a well-defined notion of what it means to `vary $\lambda$' with $z$ held fixed.  To do this, let $\lambda_1 \in \Lambda_0$, $m \ge 0$ and let $Z$ be any relatively compact open subset of $\mathcal{A}_{\lambda_1; \infty,m}$ which for now we shall assume is bounded. We claim there exists a positive real number $r$ such that $Z\subset \mathcal{A}_{\lambda;\infty,m}$ for all $\lambda \in {\mathrm B}(\lambda_1,r)$. Let $g , G$ be the minimum and maximum values respectively of $G_{\lambda_1;m}(z)$ taken on $\overline{Z}$
and let $i$ be a non-negative integer such that for all $z \in \overline{Z}$,
\vspace{-.2cm}

$$
|Q_{\lambda_1;m,m+i}(z)|>2R_{0}.
$$

\vspace{.1cm}
If we apply (3) to the formula for Green's function given earlier, we see that, for all $\lambda \in \Lambda_0$ and all $R \ge R_0$, 

\vspace{-.35cm}
\begin{equation}
\log R - (1/2)\log 2 \leq G_{\lambda;m}(z) \leq \log R + (1/2)\log(3/2), \qquad |z| = R.
\end{equation}

It is clear from the functional properties of $G_{\lambda_1;m}$ that $i$ will depend only on $g$.  Moreover, using the boundedness of $Z$ and the bound $B$ on the norms of the Fr\'echet differentials for the coefficients of the polynomials $\{P_{\lambda;n}\}_{n=m}^\infty$, $\lambda \in \Lambda_0$, we may find $r = r(g, G, B, \epsilon, ||\lambda_1 - \lambda_0||, d, M) \le \epsilon_0 - ||\lambda_1 - \lambda_0||$ such that if we set $\Lambda_1 = {\mathrm B}(\lambda_1, r)$, then for all $\lambda \in \Lambda_1$ and $z \in \overline{Z}$
\vspace{-.2cm}

$$
|Q_{\lambda;m,m+i}(z)-Q_{\lambda_1;m,m+i}(z)|<R_{0}.
$$

\vspace{.1cm}
Thus $|Q_{\lambda;m,m+i}(z)|>R_0$ and as a result, we see that all $z \in {\overline Z}$ escape under iteration by $\pemlam$ for any given $\lambda \in \Lambda_1$, and thereby must belong to the corresponding basin at infinity $\mathcal{A}_{\lambda;\infty,m}$. This tells us that we may vary $\lambda$ within $\Lambda_1$ and $Z$ will remain in the domain of all of the functions $\phimn$.

Each $Q_{\lambda;m,n}$ is holomorphic as a function of $\lambda$ on $\Lambda_0$ with $z$ held fixed since it is a finite composition of polynomials in $z$, each of which has coefficients which are holomorphic in $\lambda$. Using (3) and the maximum modulus theorem, these functions are uniformly bounded on $Z \times \Lambda_1$. For $z$ large, real and positive, for fixed $m, n$, by definition all these functions $\phi_{\lambda;m}^n$ use the same locally defined standard branch of $w^{1/D_{m,n}}$ for every $\lambda \in \Lambda_1$. For other values of $z$ in $\mathcal{A}_{\lambda;\infty,m}$, the value of $\phi_{\lambda;m}^n(z)$ is obtained by analytic continuation of the root $w^{1/D_{m,n}}$ from the positive real axis. This analytic continuation of the root is in turn determined by continuation of the standard branch of the argument $\arg Q_{\lambda;m,n}(z)$ from the positive real axis and for nearby values of $\lambda$, the values of this argument will be close. It then follows that for any fixed $z \in Z$, and any $\lambda \in \Lambda_1$, we can find a parameter neighbourhood of $\lambda$ in $\Lambda_1$ on which we may use the same locally defined branches of the root $w^{1/D_{m,n}}$ to form $\phi_{\lambda;m}^n(z)$ for each different value of $z$. From this it follows easily that the functions $\phi_{\lambda;m}^n$ will also be holomorphic in $\lambda$.

From above, each $\phi_{\lambda; m}^n(z)$ is conformal as a function of $z$ with $\lambda$ held fixed.  Since $R_{ 0}$ is an escape radius valid for all $\pemlam$ with $\lambda \in \Lambda_1 \subset \Lambda_0$, if $|z| = R_{ 0}$ then our estimate for $R_0$ ensures that for all $\lambda$ in this neighbourhood, by (3) 
\vspace{-.4cm}

\begin{equation}
 |\phi_{\lambda; m}^n(z)|  \le \sqrt{\frac{3}{2}}R_0.
\end{equation}

Repeating the argument with $cR_0$ for an appropriate constant $c>1$ and applying the maximum modulus theorem shows that $\phi_{\lambda;m}^n$ must also be uniformly bounded in $n$ and $\lambda$ on $Z \times \Lambda_1$ and that this bound must be uniform with respect to $d$, $M$ and $G$. 

If we now let $\lambda \in \Lambda_1$ and look at the ball ${\mathrm B}(\lambda, r - ||\lambda - \lambda_1||) \subset \Lambda_1$, then, as in the proof of Theorem 2.5, an application of the Schwarz Lemma to complex lines on the ball ${\mathrm B}(\lambda, (r - ||\lambda - \lambda_1||)/2)$ of half the radius implies that the G\^ateaux- and hence the Fr\'echet-derivatives of these functions are bounded for $\lambda \in \Lambda_1$ and $z \in Z$ and these bounds are uniform with respect to $\epsilon$, $B$, $d$, $M$, $\tfrac{||\lambda - \lambda_1||}{r}$, $g$ and $G$. In particular they do not depend on $n$. Thus, if we can show uniform convergence in $n$, Lemma 2.1 from \cite{Com3} gives us that the limit function will also be holomorphic in $\lambda$ and the norm of its differential will be bounded with these same bounds.

Now, for each $n > m+1$, we have the functional equation\vspace{.2cm}
\begin{eqnarray}
\phi_{\lambda ;m+1}^n \circ P_{\lambda;m+1}(z) & = & (Q_{\lambda; m+1,n}\circ P_{\lambda;m+1}(z))^{1/D_{m+1,n}}\\ 
\ & = & \left[(Q_{\lambda; m,n}(z))^{1/D_{m,n}}\right]^{d_{m+1}}\notag \\
\ & = & (\phi_{\lambda ;m}^n)^{d_{m+1}}.\notag
\end{eqnarray}

Note that all branches of roots taken in the above are the continuation from some point far out on the positive real axis of the standard ones as above.
Suppose $z \in Z$, $\lambda \in \Lambda_1$ and let $n_0$ be such that $|Q_{\lambda;m,n_0}(z)|=R \ge R_0$.  Then for $n \ge n_0$, 

\begin{eqnarray*}
\left| \frac{\phi_{\lambda;m}^{n+1}(z)}{\phimn} - 1 \right| & = & \left| \left(\frac{Q_{\lambda;m,n+1}(z)}{Q_{\lambda;m,n}(z)^{d_{n+1}}}\right)^{1/D_{m;n+1}}-1\right|\\
\ & = & \left| \left(\frac{P_{\lambda;n+1}\circ Q_{\lambda;m,n}(z)-Q_{\lambda;m,n}(z)^{d_{n+1}}}{Q_{\lambda;m,n}(z)^{d_{n+1}}} + 1 \right)^{1/D_{m,n+1}}-1\right| .
\end{eqnarray*}

Again all branches of roots taken in the above are the continuation from some point far out on the positive real axis of the standard ones as above. Using the observation that when $|u| \leq 1/2$ and $p \ge 1$,
\vspace{-.2cm}

$$
|(1+u)^{1/p}-1|<\frac{1}{p},
$$

we see that if we express $\phimn = \phi_{\lambda;m}^m \prod_{i=m+1}^n {(\phi_{\lambda;m}^i/\phi_{\lambda;m}^{i-1})}$ (where $\phi_{\lambda;m}^m$ is simply the identity), these functions $\phimn$ will converge if 

$$
\left|\frac{P_{\lambda;n+1}\circ Q_{\lambda;m,n}(z)-Q_{\lambda;m,n}(z)^{d_{n+1}}}{Q_{\lambda;m,n}(z)^{d_{n+1}}}\right| \leq \frac{1}{2}.
$$

\vspace{.2cm}
However, if $n \ge n_0$, $|Q_{\lambda;m,n}(z)| \ge R_0$ and, applying (2) to the lefthand side of this inequality, we see that this is indeed true and the desired convergence follows.

These convergence estimates for $\phimn$ are uniform with respect to the integer $n_0$ which in turn depends on the minimum $g$ of the Green's function used earlier, together with our choice of the radius $r$ of the parameter ball $\Lambda_1$. These convergence
estimates for $\phimn$ must thus hold uniformly on $Z\times \Lambda_1$, and so the limit $\phi_{\lambda;m}$ must be holomorphic in $\lambda \in \Lambda_1$ with $z \in Z$ held fixed. 

From above, the operator norms of the Fr\'echet differentials of the functions $\phimn$ as functions of $\lambda$ with $z$ held fixed are bounded on $Z \times \Lambda_1$. By Lemma 2.1 of \cite{Com3} again, it then follows that the same is true for the operator norms of the Fr\'echet differentials of the B\"ottcher maps $\phi_{\lambda;m}$. Lemma 2.1 of \cite{Com3} once more and the bounds we obtained for the Fr\'echet differentials of the functions $\phimn$ show that these bounds are then
uniform with respect to $\epsilon$, $B$, $d$, $M$, $\tfrac{||\lambda - \lambda_1||}{r}$, $g$ and $G$.

By uniform convergence and (6), the functions $\phi_{\lambda;m}$ must clearly also satisfy the functional equations $\phi_{\lambda;m+1}(P_{m+1}(z)) = (\phi_{\lambda;m})^{d_{m+1}}$ on $Z \times \Lambda_1$ for each $m \ge 0$ and a similar argument to the above shows that we can easily extend this to hold for every $\lambda \in \Lambda_0$ and every $z \in \mathcal{A}_{\lambda;\infty,m}$. 

We still need to show that we can remove the assumption that $Z$ is bounded and also show that the bounds on the operator norms of the Fr\'echet derivatives of the functions 
$\phi_{\lambda;m}$ as functions of $\lambda$ do not depend on the maximum value of the Green's function $G_{\lambda_1; m}$. Since all the polynomials in our analytic family are monic, for each $\lambda \in \Lambda_0$, $\phi_{\lambda;m}$ must have $z$-derivative $1$ at infinity and thus 
be conformal in the $z$ variable in a neighbourhood of $\infty$ with $\lambda$ held fixed.  Utilising a lifting argument to pull back via the functional equation then allows us to conclude that $\phi_{\lambda;m}$ will be conformal in $z$ and map $\mathcal{A}_{\lambda;\infty,m}$ to $\sphere \setminus \overline \D$. We remark that from this it is then straightforward to check that $\log|\phi_{\lambda;m}|$ agrees with the Green's function $G_{\lambda;m}(z)$ with pole at infinity on the domain $\mathcal{A}_{\lambda;\infty,m}$.

Recall the classical class $\mathcal S$ of all univalent functions $f$ on the unit disc which satisfy $f(0) = 0$, $f'(0) = 1$. By Theorem 1.8 on page 4 of \cite{CG}, we have a weak form of the Bieberbach conjecture which states that if 
$f(z) = z + \sum_{n=2}^\infty a_n z^n \in {\mathcal S}$, then $|a_n| < en^2$ for each $n \ge 2$.  All the functions $\phi_{\lambda;m}$ are univalent on $\sphere \setminus \overline {\mathrm D}(0, R_0)$ and have derivative $1$ at infinity. It then follows easily that we can find $D > 0$ such that on $\sphere \setminus \overline {\mathrm D}(0, 2R_0)$, we have $|\phi_{\lambda; m}(z) - z| < D$. This shows that the Fr\'echet derivatives of the functions $\phi_{\lambda;m}$ remain bounded near infinity and that, for $|z| > 2R_0$, these bounds will depend only on $\tfrac{||\lambda - \lambda_0||}{\epsilon_0}$ (and in particular not on $r$ as $\sphere \setminus \overline {\mathrm D}(0, R_0)$ is a common domain for all the functions $\phi_m$) and the degree and coefficient bounds $d$ and $M$ (via the escape radius $R_0$). 

Finally, note that, in the case of a monic analytic family, if we had another sequence 
of B\"ottcher mappings $\tilde \phi_{\lambda;m}$ defined on ${\mathcal A}_{\lambda; \infty ,m}$ for each $\lambda \in \Lambda_0$, then the compositions ${\tilde \phi_{\lambda;m}}^{\circ -1} \circ \phi_{\lambda;m}$ would give us a family of conformal mappings of $\sphere \setminus \overline \D$ to itself all of which had derivative $1$ at infinity. This implies that all these compositions must be the identity from which the uniqueness of the B\"ottcher mappings in the monic case follows. 
With this, we have (i) and (iii). 

Let $m \ge 0$ and let $\psi_{\lambda;m}$ be the $z$-inverse of $\phi_{\lambda;m}$ so that $\phi_{\lambda;m}(\psi_{\lambda;m}(z))=z$ for all $z \in \sphere \setminus \overline{\mathbb{D}}$.  From above these mappings are clearly uniquely defined for a monic analytic family. We now verify that this $z$-inverse is Fr\'{e}chet-holomorphic in $\lambda$ as well. Let $W$ be a relatively compact bounded open subset of $\sphere \setminus \overline \D$, let $\lambda_1 \in \Lambda_0$ and let $Z \subset \mathcal{A}_{\lambda_1;\infty,m}$ be the image of $W$ under $\psi_{\lambda_1;m}$. Next, let $r > 0$ and let $\Lambda_1 = {\mathrm B}(\lambda_1, r) \subset \Lambda_0$ be defined similarly to above so that all the functions $\phi_{\lambda; m}$, $\lambda \in \Lambda_1$ are defined on $Z$.

Let $\mathbf{\xi}$ be a fixed but arbitrary unit vector in $\ell^\infty(\mathbb{C})$ and define the function $f_m:(Z \times \mathbb{D}) \rightarrow (\sphere\setminus \overline{\mathbb{D}})$ by 
\vspace{-.5cm}

$$
f_m:(z,t)\mapsto \phi_{\lambda_{1} + r  t \mathbf{\xi};m}(z).
$$

For $z_0 \in Z$, we have $\partial f_m/\partial z \neq 0$ at $(z_0,0)$.  If $f_m(z_0,0)=w_0$, then, by the Implicit Function Theorem, there exists an open bidisk $U \times V$ containing $(z_0,0)$ and a unique analytic function $u_m(t)$ such that $f_m(z,t)=w_0$ on $U\times V$ if and only if $z=u_m(t)$ for some $t \in V$.  By definition, we also have $f(\psi_{\lambda_1+\epsilon t\mathbf{\xi} ;m}(w_0),t)=w_0$ for all $t \in \mathbb{D}$.

By the uniqueness of the implicit function $u_m(t)$ guaranteed by the Implicit Function Theorem, $u_m(t)$ and $\psi_{\lambda_{1}+\epsilon t\mathbf{\xi} ;m}(z_0)$ agree on $V$.  In particular $\psi_{\lambda_{1} + \epsilon t \mathbf{\xi} ;m}(w)$ is analytic with respect to $t$ with $w$ held fixed, and so $\psi_{\lambda ;m}(w)$ is G\^ateaux-holomorphic with respect to $\lambda$ at each $w \in \sphere\setminus \overline{\mathbb{D}}$.  

$\overline{W}$ is compact, and so we may find $R > 1$ so that $W\subset {\mathrm D}(0,R)$.  By the existence of the uniform escape radius $R_0$ for every sequence $\pemlam$, $\lambda \in \Lambda_0$ and the fact that the functions $\psi_{\lambda;m}$ are isometries of the hyperbolic metric, it follows that they will be uniformly bounded in $\lambda$ and $m$ on $\overline W$. By replacing $\epsilon_0$ by $\epsilon_0/2$ in the definition of $\Lambda_0$, using the above local boundedness in $\lambda$ and applying an argument using the Schwarz Lemma similar to the case for $\phi_{\lambda;m}$ we see that the Fr\'echet differential $\mathbf{d}_\lambda\psi_{\lambda;m}$ is a bounded operator and its norm is bounded uniformly in terms of $d$, $M$, $\tfrac{||\lambda - \lambda_0||}{\epsilon_0}$, $B$ and the maximum value of $\log |z|$ on $W\times\Lambda_0$. 

Note, however, that the bound on the absolute value of $\psi_{\lambda;m}(z)$ does not depend on how close $z$ is to $\overline \D$ and so the bound on the operator norm of the Fr\'echet differential will not depend on the minimum value of $\log |z|$ on $W$.
Looking at the other extreme, a very similar argument as before shows that we can remove the assumption that $W$ is bounded and that the Fr\'echet differentials $\mathbf{d}_\lambda\psi_{\lambda;m}$ are uniformly bounded near infinity and that these bounds will depend only on $\tfrac{||\lambda - \lambda_0||}{\epsilon_0}$ and $d$ and $M$. Thus, we may as well assume that $W$ is all of $\sphere \setminus \overline \D$ and with this, we have (ii) and (iv) and the proof is now complete. $\blacksquare$

We remark that, using Riemann maps and normal families, Sester in \cite{Ses2} Proposition 2.3 earlier proved a result on continuity of B\"ottcher maps for quadratic fibered polynomials which overlaps somewhat with our setting. However, for analytic rather than continuous dependence, more work is required. 

From the above, we have the following corollary that says the external rays move holomorphically in the appropriate sense. Similarly to the classical case, we may define an external ray of angle $\theta$ at time $m$

\vspace{-.25cm}
$$
R_{\theta, m}=\{\psi_m(Re^{i\theta}):R>1\}
$$

\vspace{0cm} 

where $\psi_m$ are the inverse B\"ottcher mappings introduced in Theorem 3.1. For an analytic family ${\mathcal P}_\Lambda$ as above, if $\lambda \in \Lambda$ and $m \ge 0$, we use the notation $R_{\lambda; \theta, m}$ to denote the corresponding external ray with angle $\theta$ at time $m$ for the sequence $\pemlam$. 

\vspace{.2cm}
\begin{corollary}
Let $\Lambda$, ${\mathcal P}_\Lambda$, $\lambda_0$ be as in the previous result and let $\Lambda_0$ be the ball of radius $\epsilon_0$ about $\lambda_0$ obtained as above. Then, for any $m \ge 0$, the iterated basin at infinity ${\mathcal A}_{\lambda_0 ; m}$ and any external ray $R_{\lambda_0; \theta, m}$, $0 \le \theta < 2\pi$, move holomorphically on $\Lambda_0$ and the Fr\'echet derivative of this holomorphic motion depends only on $d$, $K$, $M$, $\tfrac{||\lambda - \lambda_0||}{\epsilon}$ as above and in particular not on the point chosen on $R_{\lambda_0; \theta, m}$.

Furthermore, if ${\mathcal P}_\Lambda$ is monic or globally analytically conjugate to a monic analytic family and $\lambda$ belongs to a hyperbolic component $\Omega$, then these holomorphic motions can be extended to give us a holomorphic motion on all of $\Omega$. If ${\mathcal P}_\Lambda$ is a monic family, then this holomorphic motion is uniquely defined.
\end{corollary}

\textbf{Proof:} The proof of the first part is immediate from Theorem 3.1 and the fact that for $\lambda \in \Lambda_0$ we have ${\mathcal A}_{\lambda; m} = \psi_{\lambda; m} \circ \phi_{\lambda_0; m} ({\mathcal A}_{\lambda_0 ; m})$ and $R_{\lambda; \theta, m} = \psi_{\lambda; m} \circ \phi_{\lambda_0; m} (R_{\lambda_0; \theta, m})$ (note that using the standard estimates for the hyperbolic metric, e.g. Theorem 4.3 on P. 13 of \cite{CG}, it is not hard to find three points $x_1$, $x_2$, $x_3$ on $R_{\lambda; \theta, m}$ for which condition (iv) for a holomorphic motion is fulfilled). 

For the second part of the statement, we first assume that  ${\mathcal P}_\Lambda$ is a monic family. The 
B\"ottcher maps $\phi_{\lambda; m}$, $m \ge 0$ and hence their inverses $\psi_{\lambda; m}$, $m \ge 0$ are then uniquely defined on all of $\Omega$ by Theorem 3.1. These give a locally defined holomorphic motion which by a similar argument to the proof of Proposition 3.1 can be continued along any path from $\lambda_0$ to any other point of $\Omega$. Since all the mappings in this holomorphic motion are analytic in $z$ with $z$-derivatives $1$ at infinity, again by Theorem 3.1, it follows that this motion is uniquely defined on all of $\Omega$. For a general analytic family, the result follows by applying the global analytic conjugacy to a monic family. $\blacksquare$

Before we close this section, we remark that although Theorem 3.1 seemed to have nothing to do with hyperbolicity and the assumption that all the iterated Julia sets were connected for our family ${\mathcal P}_\Lambda$ seemed to be a long way from saying that the corresponding sequences were hyperbolic, this is not in fact the case. 

\begin{theorem}
Let ${\mathcal P}_\Lambda$ be an analytic family over an open subset $\Lambda$ of $\linf$ and let $\lambda$ be in the interior of ${\mathcal M}_\Lambda$. Then $\lambda \in \mathcal{HC}_\Lambda$. 
\end{theorem}

\textbf{Proof:} If the sequence $\pemlam$ is not hyperbolic, then, by Theorem 2.3, the postcritical distance is zero. This implies that we can make arbitrarily small perturbations where we change just one polynomial of the sequence slightly so as to move some of the critical values into the basin of infinity at some time $m \ge 0$. However, in view of our assumption that $\lambda$ is in the interior of the connectedness locus, this is clearly impossible. $\blacksquare$

Effectively, the above states that there can be no `queer components' in non-autonomous iteration, at least in the connected case. Note that this does \emph{not} imply that there can be no components of the connectedness locus in the classical case where we have a persistently non-hyperbolic point as described in \cite{MSS}. The basic reason is that there is far more freedom available for making a non-autonomous perturbation than there is in the classical case. Note also that we are also not stating that hyperbolicity is dense in the connectedness locus. There could still possibly be `empty filaments' of ${\mathcal M}_\Lambda$ which are far away from points of $\mathcal{HC}_\Lambda$.

We finish with a lemma on the lengths of segments of hyperbolic geodesics, which in our case are external rays. This is a version of a result in \cite{CJY} (Lemma 3.3) which is used to show the basin of infinity is a John domain for a semi-hyperbolic polynomial in the classical case. For a non-autonomous version, a proof of this can be found in the work of Sumi, specifically Claim 4 in the proof of Theorem 1.12 in \cite{Sum2} (note that this also requires the estimate on the Green's function in terms of the distance to the boundary which is proved as claim ($*$) on the second page of \cite{Sum3}). For a point $z$ in the basin at infinity at time $m \ge 0$ for a bounded sequence with connected iterated Julia sets, following Carleson, Jones and Yoccoz, let us denote by $\gamma_z$ the segment of the Green's line in ${\mathcal A}_{\infty, m}$ which runs from $z$ to $\partial {\mathcal A}_{\infty, m} = {\mathcal J}_m$ (and which in our case will be part of an external ray). Finally, for a curve $\gamma$, we denote the Euclidean arc length of $\gamma$ by $\ell(\gamma)$. 

\vspace{.2cm}
\begin{lemma} Let $\Pm$ be a bounded hyperbolic sequence all of whose iterated Julia sets are connected. Then there exist constants $C > 0$, $\alpha >0$ depending only on the hyperbolicity bounds for $\Pm$ such that for any $m \ge0$ and $z \in {\mathcal A}_{\infty, m}$, $\ell (\gamma_z) \le C G_m(z)^\alpha$. 
\end{lemma}

The following is an immediate consequence of this lemma.

\vspace{.2cm}
\begin{corollary}
Let $\{P_m\}_{m=1}^\infty$ be a hyperbolic bounded sequence all of whose iterated Julia sets are connected. Then for every $m \ge 0$, the iterated basin at infinity ${\mathcal A}_{\infty, m}$ is locally connected and every external ray in ${\mathcal A}_{\infty, m}$ lands at just one point on ${\mathcal J}_m$. 
\end{corollary}

Note that it follows immediately from the above local connectivity that every point on each of the iterated Julia sets for such a sequence is the landing point of at least one external ray. 
We also remark here that, using the fact as mentioned above that, for a hyperbolic bounded sequence all of whose iterated Julia sets are connected the iterated basins of infinity are John domains, one can show easily that at most finitely many external rays may land at any given point on one of the iterated Julia sets.

\section{Combinatorial Rigidity and Preservation of External Rays}

With our technical preliminary results out of the way, we are now in a position to formally state and prove our main result on external rays.  Basically, this states that the landing points of rays give us the same holomorphic motion of the iterated Julia sets on the components of $\mathcal{HC}_\Lambda$ as the holomorphic motions of Theorem 2.4 which were constructed using grand orbits. Further, if a number of external rays meet at a given point, then we have combinatorial rigidity in the sense that this picture is preserved on the entire hyperbolic component, including the angles of the rays. 

Before stating our result, we remark that Sester in Section 2 of \cite{Ses2} proves a local combinatorial rigidity result about the continuous dependence of the landing point of 
external rays for quadratic fibered polynomials. In our case, we show analytic dependence not just of the landing point, but whole external rays. Further, our result is valid on entire hyperbolic components and the necessary estimates on the Green's function are a consequence of hyperbolicity via Lemma 3.2 rather than assumptions as in Proposition 2.7 of \cite{Ses2}.

\vspace{.2cm}
\begin{theorem}  Let $\Lambda \subset \linf $ be open and let ${\mathcal P}_{\Lambda}$ be an analytic family of bounded polynomial sequences over $\Lambda$ which is either monic or globally analytically conjugate to a monic analytic family. Let $\lambda_0 \in \mathcal{HC}_\Lambda$ with $\lambda_0$ belonging to some hyperbolic component $\Omega \subset \Lambda$ in parameter space. We then have the following.

\begin{enumerate}

\item The holomorphic motion of the iterated basins of infinity ${\mathcal A}_{\lambda; m}$, $\lambda \in \Omega$ of Corollary 3.1 extends via the landing points of external rays to a holomorphic motion of the sets ${\mathcal A}_{\lambda; m} \cup {\mathcal J}_{\lambda; m}$ and its restriction to the iterated Julia sets ${\mathcal J}_{\lambda; m}$ coincides with a uniquely defined extension to all of $\Omega$ of the locally defined holomorphic motion $\tau_{\lambda; m}(z)$ of Theorem 2.4 which is initially defined in a neighbourhood of $\lambda_0$. In particular $\tau_{\lambda; m}(z)$ can be extended in a unique fashion from a neighbourhood of $\lambda_0$ to all of $\Omega$.  

\vspace{.2cm}
\item Let $p \in {\mathcal{J}_{\lambda_0;0}}$ be such that the $N$ external rays $R_{\lambda_0;\theta^{1,m}, m}, R_{\lambda_0;\theta^{2,m}, m},..., R_{\lambda_0;\theta^{N,m}, m}$ (and no others) meet at $p_m=Q_{\lambda_0;m}(p)$ for each $m \ge 0$.  
Then, for each $m \ge 0$, $\lambda \in \Omega$, the $N$ external rays $R_{\lambda;\theta^{1,m}, m}, R_{\lambda;\theta^{2,m}, m},..., R_{\lambda;\theta^{N,m}, m}$ with these same angles and no other rays meet at $\tau_{\lambda; m}(p_m)$.

\end{enumerate}

\end{theorem}

\textbf{Proof:}  We start by noting that in the monic case, in view of the uniqueness part of Theorem 3.1, we can assume that the B\"ottcher maps $\phi_{\lambda; m}$ and their inverses can then be uniquely defined on all of $\Omega$. For a general analytic family, we may use the global analytic conjugacy to a monic sequence to ensure that we have B\"ottcher maps $\phi_{\lambda; m}$ which give us well-defined functions in $\lambda$ on all of $\Omega$. In particular, external rays $R_{\lambda; \theta, m}$ also give us well-defined functions of $\lambda$. 

Let us therefore assume from now on that ${\mathcal P}_{\Lambda}$ is monic and let $\lambda_1 \in \Omega$.  $\Omega$ is connected and so we may find a continuous path $\gamma:[0,1]\rightarrow \Omega$ with $\gamma(0)=\lambda_0$ and $\gamma(1)=\lambda_1$.  Again for each $t \in [0,1]$, $\lambda_t:=\gamma(t)$ represents a hyperbolic sequence, and as such by Theorem 2.1 we may find a ball ${\mathrm B}(\lambda_t,\epsilon_t) \subset \Omega$ so that each parameter in this ball gives rise to a polynomial sequence and these sequences are all uniformly bounded and hyperbolic with the same constants.  As in the proof of Proposition 3.1 we may then take a finite subset of these balls to cover $\gamma$ so that we have an open parameter neighbourhood $X \subset \Omega$ of $\gamma$ of uniformly bounded and hyperbolic sequences (with the same constants).

Since by Theorem 2.3 the postcritical distance depends only on the hyperbolicity bounds, this further guarantees that there exists a universal lower bound $\delta >0$ for the postcritical distance on $X$. 
By making the discs of $X$ smaller if needed (in which case we may need more, but still finitely many of them) so that we may apply Theorem 2.4 on each one, by moving along successive balls, we may continue the embeddings $\tau_{\lambda; m}(z)$ along $\gamma$ (note that, by either using the grand orbit of the continuation of the same point in the iterated Julia sets or by Theorem 2.5, we can ensure that this continuation will agree on the overlap between any two successive balls). As these are homeomorphisms in $z$ and all the iterated Julia sets ${\mathcal J}_{\lambda_0; m}$, $m \ge 0$ for $\lambda_0$ are connected, all the iterated Julia sets ${\mathcal J}_{\lambda;m}$, $m \ge 0$ are also connected for any $\lambda \in X$.

Now let $x^1 = 2$ and for each $i > 1$, let the points $x^i$ be chosen such that $|x^i| < |x^{i-1}|$ with $\rho_{\sphere \setminus \overline \D} (x^i, x^{i-1}) = 1$. For $0 \le \theta < 2\pi$ fixed and $i \ge 1$, let $z^i$ be the point $x^i e^{\i \theta}$. Using this, for each $\lambda \in X$, and each $m \ge 0$, $i \ge 1$, $0 \le \theta < 2\pi$ define $p^i_m(\lambda)$ to be the point 
$\psi_{\lambda; m}(z^i)$. By (iv) of Theorem 3.1 and Lemma 3.2, we can again make the balls of $X$ smaller if needed so that for each $m \ge 0$ the functions $p^i_m(\lambda)$ converge uniformly on $X$ as $i \to \infty$ to a limit function $p_m(\lambda)$. Again by (iv) of Theorem 3.1 and Lemma 2.1 of \cite{Com3}, these limit functions will be Fr\'echet-holomorphic on $X$ and the operator norms of their Fr\'echet differentials will be bounded uniformly in terms of the uniform degree and coefficient bounds for  the sequences on $X$. 

By Corollary 3.2 the functions $p_m(\lambda)$ are the landing points of the external rays $R_{\lambda; \theta, m}$ and since these are well-defined on $\Omega$ from above, this allows us to define a Fr\'echet-holomorphic function on all of $\Omega$. If, for $\lambda_0$, $R_{\lambda_0; \theta, m}$ lands at $p_m$, then 
by making the balls of $X$ smaller yet again if needed so that we can apply Theorem 2.5 on each of them, by applying this result on successive balls and using the fact that the functions $p_m(\lambda)$ are well-defined on $\Omega$, we have that $p_m(\lambda) = \tau_{\lambda; m}(p_m)$ on $\Omega$ where $\tau_{\lambda;m}$ is the continuation of the holomorphic motion of Theorem 2.4 from a neighbourhood of $\lambda_0$ to all of the set $X$ as constructed above (note that, because the functions $p_m(\lambda)$ are well-defined, in particular this will not depend on our choice of set $X$ for joining $\lambda_0$ to any other point of $\Omega$). Theorem 2.5 shows that this is unique for a given path from $\lambda_0$ to $\lambda_1$ and the fact that the functions $p_m(\lambda)$ are well-defined on $\Omega$ shows that this extension is uniquely defined on $\Omega$ which proves the first part of the statement. The second part then follows immediately from this. $\blacksquare$ 

Lastly, we have the following corollary.

\begin{corollary}
Let $\Lambda \subset \linf$ be open and let ${\mathcal P}_{\Lambda}$ be an analytic family of bounded polynomial sequences over $\Lambda$.

Let $\Omega$ be a hyperbolic component for ${\mathcal P}_{\Lambda}$, let $\lambda_0 \in \Omega$ and suppose that we can find 
a holomorphic motion $\sigma_{\lambda; m}(z)$, $m \ge 0$ of the iterated Julia sets which gives a conjugacy on the iterated Julia sets on some neighbourhood of $\lambda_0$ in the sense that (4) and (5) of Theorem 2.4 above are satisfied. Suppose also that for any $\lambda_1 \in \Omega$, and any path $\gamma$ in $\Omega$ from $\lambda_0$ to $\lambda_1$, we may continue $\sigma_{\lambda; m}(z)$, $m \ge 0$ as a conjugacy on iterated Julia sets along a suitable neighbourhood of $\gamma$ in $\Omega$. 

Suppose further that ${\mathcal P}_{\Lambda}$ is either monic or globally analytically conjugate to a monic analytic family and that in addition at least one of the following holds: (I) $\Omega$ is a component of $ \mathcal{HC}_\Lambda$ or (II) $\Omega$ is simply connected. Then $\sigma_{\lambda; m}(z)$ can be defined on all of $\Omega$ and this extension must be unique. Also, any two such locally defined holomorphic motions must agree on all of $\Omega$. 
\end{corollary}

\textbf{Proof:}  The first case is covered by an argument similar to that at the end of the proof of the last result. The second is covered by an argument along the lines of the monodromy theorem similar to that in Proposition 3.1. Note that in view of Theorem 2.5, the holomorphic motions are locally unique on a neighbourhood of $\lambda_0$. We can then use Theorem 2.5 again or the identity principle as before so we can conclude the continuation along any fixed path is unique while the argument using homotopy along the lines of the monodromy theorem shows that this continuation is independent of the choice of path. $\blacksquare$

An interesting contrast to Theorem 4.1 above can be seen in the paper of Laura DeMarco and Suzanne Lynch Hruska \cite{DeMH}. In Section 5 of their paper they consider the family $F_a(z,w) = (z^2, w^2 + az)$ of polynomial skew-products and of particular interest is their picture for the case $a= -1$, as illustrated in Figure 1 of that paper. For $z = 1$ we have the classical `basilica' and this rotates through an angle of $\pi$ as $z$ moves around the unit circle, which clearly does not preserve the angles of external rays. However, this is not an analytic family in the non-autonomous sense as the angle doubling properties of $z^2$ on the unit circle show that the $z$-derivatives of the constant coefficients quickly become unbounded. To be precise, if we pick some fixed $z \in {\mathrm C}(0,1)$ and let $P_1(w) = w^2 + az$, then $P_m(w) = w^2 + az^{2^{m-1}}$ for $m \ge 1$. Finally, we note in passing that DeMarco and Hruska do also prove a result on holomorphic motions for polynomial skew-products which depend analytically on a parameter (\cite{DeMH} Theorem 4.2).

While the lack of periodic orbits in the present setting prevents a more literal analogy to the standard definition of an orbit portrait, this result demonstrates that classical orbit portraits still have some bearing here.  We immediately see, for example, that if a hyperbolic component in non-autonomous parameter space meets a constant sequence $\{P,P,...\}$, then each sequence in this component possesses the same ray landings guaranteed by the orbit portrait of $P$ even though the landing points no longer constitute a recurrent orbit.  This suggests the existence of a more general structure in the non-autonomous parameter space which describes `nearly periodic' behavior and is also preserved by deformation through hyperbolic maps, but which does not rely on the existence of periodic orbits.  

\section{Questions for Further Study}

Although we have not yet established a proof, the authors believe that Theorem 4.1 can be strengthened. It seems likely that the common landing point of external rays continues to move holomorphically on a larger set than just the hyperbolic component as in the classical case for quadratic polynomials where it is a repelling periodic point on the wake for a given orbit portrait (\cite{Mi2} Theorem 1.2). 

Also of interest is whether new, non-classical orbit `portraits' can arise from a non-autonomous sequence.  It is fairly straightforward to show that, even in the non-autonomous setting, all sets of multiple rays landing at a common point are finite and are at least conformally conjugate to a set of rays with rational angles.  This seems to hint that there are no non-classical orbit portraits, although the combinatorics of the ray landings are very complicated in the case where the degrees of the polynomials in the sequence are allowed to vary.

\end{document}